\documentclass[12pt, oneside, reqno]{amsart}

\usepackage{amssymb}
\usepackage{amsmath}
\usepackage{pst-all}
\usepackage{geometry}

\evensidemargin\oddsidemargin

\begin{document}

\setcounter{page}{1}

\newtheorem{PROP}{Proposition}
\newtheorem{REMS}{Remark}
\newtheorem{LEM}{Lemma}
\newtheorem{THE}{Theorem}
\newtheorem{THEA}{Theorem A\!\!}
\newtheorem{THEB}{Theorem B\!\!}

\renewcommand{\theTHEA}{}
\renewcommand{\theTHEB}{}

\newcommand{\eqnsection}{
\renewcommand{\theequation}{\thesection.\arabic{equation}}
    \makeatletter
    \csname  @addtoreset\endcsname{equation}{section}
    \makeatother}
\eqnsection

\def\a{\alpha}
\def\bt{\beta}
\def\B{{\bf B}} 
\def\cC{{\mathcal{C}}} 
\def\CC{{\mathbb{C}}} 
\def\cL{{\mathcal{L}}} 
\def\dab{d_{\a, \bt}}
\def\Dab{D_{\a, \bt}}
\def\Hab{H_{\a, \bt}}
\def\hab{h_{\a, \bt}}
\def\mab{\mu_{\a, \bt}}
\def\cDab{{\overline \Dab}}
\def\Ea{E_\a}
\def\Eab{E_{\a, \bt}}
\def\Fab{F_{\a, \bt}}
\def\EE{{\mathbb{E}}} 
\def\fa{f_\a}
\def\tfa{{\widetilde \fa}}
\def\Fa{F_\a}
\def\gaa{g_\a}
\def\tga{{\widetilde \gaa}}
\def\ga{\gamma}
\def\Ga{\Gamma}
\def\G{{\bf \Ga}} 
\def\ha{h_\a}
\def\tha{{\widetilde \ha}}
\def\thab{{\widetilde \hab}}
\def\Ha{H_\a}
\def\tHa{{\widetilde \Ha}}
\def\i{{\rm i}}
\def\I{{\bf I}}
\def\K{{\bf K}}
\def\ka{k_\a}
\def\tka{\widetilde \ka}
\def\Ka{{\bf K}_\a}
\def\L{{\bf L}}
\def\Lb{{\mathcal L}_\bt}
\def\lbd{\lambda}
\def\lcr{\left[}
\def\lpa{\left(}
\def\lva{\left|}
\def\M{{\mathcal M}}
\def\Ma{{\mathcal M}_\a}
\def\Mla{{\bf ML}_\a}
\def\MMa{{\bf M}_\a}
\def\CMla{\widetilde \Mla}
\def\cMMa{\widetilde \MMa}
\def\cMa{\widetilde \Ma}
\def\na{n_\a}
\def\pb{{\mathbb{P}}}
\def\rl{{\mathbb{R}}}
\def\rpa{\right)}
\def\rcr{\right]}
\def\rva{\right|}
\def\S{{\bf S}}
\def\Ta{{\bf T}_\a}
\def\X{{\bf X}}
\def\tX{\widetilde X}
\def\Xa{{\bf X}_\a}
\def\Ya{{\bf Y}_\a}
\def\Sa{S_\a}
\def\U{{\bf U}}
\def\Ua{{\bf U}_\a}
\def\Un{{\bf 1}}
\def\Va{{\bf V}_\a}
\def\Wa{{\bf W}_\a}
\def\cZ{{\mathcal{Z}}} 
\def\Y{{\bf Y}}
\def\Z{{\bf Z}}
\def\Za{{\bf Z}_\a}
\def\cZa{\widetilde \Za}

\def\claw{\stackrel{d}{\longrightarrow}}
\def\elaw{\stackrel{d}{=}}
\def\qed{\hfill$\square$}
                  
\title[Mittag-Leffler functions and complete monotonicity]
      {Mittag-Leffler functions and complete monotonicity}
      
\author[Thomas Simon]{Thomas Simon}

\address{Laboratoire Paul Painlev\'e, Universit\'e Lille 1, Cit\'e Scientifique, F-59655 Villeneuve d'Ascq Cedex. {\em Email} : {\tt simon@math.univ-lille1.fr}}

\address{Laboratoire de physique th\'eorique et mod\`eles statistiques, Universit\'e  Paris Sud, B\^atiment 100, 15 rue Georges Cl\'emenceau, F-91405 Orsay Cedex. }



\begin{abstract} We consider two operations on the Mittag-Leffler function which cancel the exponential term in the expansion at infinity, and generate a completely monotonic function. The first one is the action of a certain differential-difference operator, and leads to a characterization via some necktie domain. The second one is the subtraction of the exponential term itself multiplied by an incomplete Gamma function. These results extend previous works by various authors.
\end{abstract}

\maketitle
 
\section{Introduction}

The classical Mittag-Leffler function is the entire function
$$\Ea (z) \; =\; \sum_{n\ge 0}\; \frac{z^n}{\Ga (1+\a n)}, \quad z \in \CC,\, \a > 0,$$
and can be viewed as an extension of the exponential function. The generalized Mittag-Leffler function writes
$$\Eab (z) \; =\; \sum_{n\ge 0} \; \frac{z^n}{\Ga (\bt+\a n)}, \quad z \in \CC, \, \a, \bt > 0.$$
Introduced for analytical purposes by Mittag-Leffler and Wiman at the beginning of the twentieth century, these functions have been the object of many studies. We refer to Chapter XVIII in \cite{E3} and Chapter 3 in \cite{D} for classical properties, and also to the survey \cite{HMS} for a more recent account. Nowadays these functions play an important r\^ole in fractional calculus, and find some applications in physics \cite{Mai}.  

It has been shown by Pollard  \cite{P} that the function 
$$x\,\mapsto\, \Ea(-x), \qquad x\in\rl^+,$$ 
is completely monotonic (CM) for any $\a \in(0, 1]$. Recall that a smooth function on $(0, +\infty)$ is CM if its successive derivatives have an alternating sign, starting positive. Bernstein's theorem - see e.g. \cite{W} p.160 - states that a function $f$ is CM if and only if it writes
$$f(x)\; =\; \int_0^\infty e^{-xt}\, \mu (dt)$$
for some positive $\sigma-$finite measure $\mu.$ Pollard's result was improved by Schneider \cite{Sch}, who showed that $E_{\a,\bt}(-x)$ is CM if and only if $\a\in (0,1]$ and $\bt\ge \a.$ A short proof of this latter result, using an Abelian transformation, has been obtained in \cite{Mi}. On the other hand, the author observed in \cite{TS0} that the function
$$x\,\mapsto\, \Ea(x^\a)\, -\, \a x^{\a -1}\Ea'(x^\a)$$ 
is CM for all $\a\in [1,2].$ In the present paper we pursue these lines of research and display further CM properties of the functions $\Ea$ and $\Eab.$ For every $\bt > 0,$ consider the following differential-difference operator 
$$\Lb f(x)\; =\; f'(x) \; +\; \frac{\bt -1}{x}\, (f(x) - f(0)).$$
For every $x > 0,$ set $\Fab(x) = \Eab (x^\a),$
$$\Dab (x)\; = \; \Lb\Fab(x)\; -\; \Fab(x),$$ 
and $\cDab (x) = -\Dab(x).$ It is easy to see that $D_{1, \bt} = 0$ for all $\bt > 0.$ If $\a \neq 1,$ notice that $\Dab$ is the difference of two functions with exponential growth at infinity: one has the convergent series representation
\begin{equation}
\label{CVSR}
\Dab (x) \; =\;  \sum_{n\ge 1} \frac{x^{\a n-1}}{\Ga (\bt+\a n-1)}\; -\;\sum_{n\ge 0} \frac{x^{\a n}}{\Ga (\bt+\a n)}\cdot
\end{equation}

\medskip

\begin{THEA} Assume $\a\neq 1.$ The following equivalences hold.

\medskip

{\em (a)} $\Dab$ {\em is CM} $\;\Leftrightarrow\;\a \in (0,1), \,\bt \ge \a\vee (1-\a).$ 

\smallskip

{\em (b)} $\cDab$ {\em is CM} $\;\Leftrightarrow\;$ {\em either} $\a \in (1,2], \, \bt\ge 1$ {\em or} $\a\in (0,1), \, \bt \le \a\wedge (1-\a).$ 

\end{THEA}

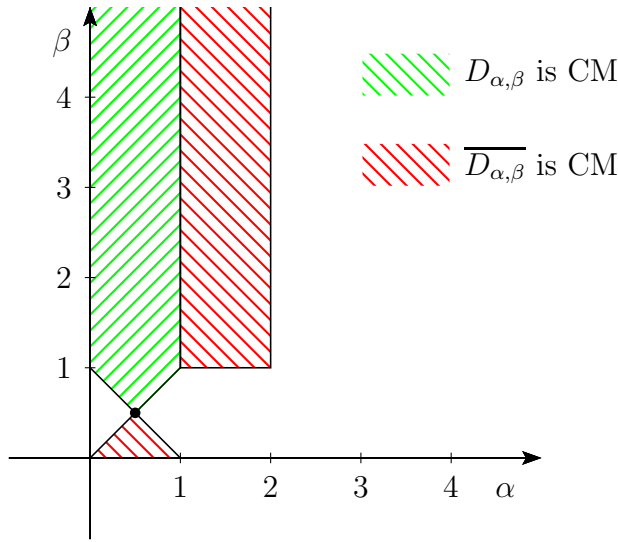
\bigskip

\begin{figure}[h]
\centering

\definecolor{vertfonce}{rgb}{0,0.5,0}

\psset{unit=1.2cm,algebraic=true,linewidth=0.6pt,
,arrowsize=4pt 2,arrowinset=0.2}

\begin{pspicture*}(-1,-1)(8,5)

\rput[t](4.6,-0.28){$\alpha$}
\rput[r](-0.2,4.6){$\beta$}

\pspolygon[hatchcolor=green,fillstyle=hlines,hatchangle=45,hatchsep=0.1](0,1)(0,10)(1,10)(1,1)(0.5,0.5)
\pspolygon[hatchcolor=red,fillstyle=hlines,hatchangle=-45,hatchsep=0.1](0.5,0.5)(1,0)(0,0)
\pspolygon[hatchcolor=red,fillstyle=hlines,hatchangle=-45,hatchsep=0.1](1,1)(1,10)(2,10)(2,1)

\psframe[linecolor=white,hatchcolor=green,fillstyle=hlines,hatchangle=-45,hatchsep=0.1](3,4)(4,4.5)
\uput[r](4,4.25){$\Dab$ is CM}
\psframe[linecolor=white,hatchcolor=red,fillstyle=hlines,hatchangle=-45,hatchsep=0.1](3,3)(4,3.5)
\uput[r](4,3.25){$\cDab$ is CM}

\psdots[dotstyle=*,dotsize=4pt](0.5,0.5)

\psaxes[Dx=1,Dy=1,ticksize=-2pt 2pt]{->}(0,0)(-0.9,-0.9)(5,5)

\end{pspicture*}

\end{figure}

\medskip

\noindent
Observe that $D_{1/2,1/2} = 0,$ as can also be seen from (\ref{CVSR}). In the above result, the r\^ole of the operator $\Lb$ is to cancel the leading exponential term in the expansion of $\Eab (x^\a).$ More precisely, the asymptotic expansion 18.1(22) p. 210 in \cite{E3} shows that for $\a\in (0,2],$
$$\Eab(x^\a)\; \sim\; \frac{x^{1-\bt}e^{x}}{\a}, \qquad x\to +\infty,$$
and the right hand side is, up to some function with polynomial decay, annihilated by the action of ${\rm Id} - \Lb.$ Other differential-difference operators can be chosen in order to make this cancellation, but $\Lb$ is the most natural one because $\Lb E_{1,\bt} = E_{1, \bt}.$ 

In view of the complete expansion 18.1(22) in \cite{E3}, one may ask if subtracting from $\Eab (x^\a)$ the leading exponential term {\em itself} would not lead to a CM function. The following result shows that this is indeed the case for $\a\in (0,2], \bt \ge 1,$ up to a slight multiplicative correction when $\bt > 1.$ For every $u,x >0,$ set
$$\ga(u,x)\; =\; \int_0^x t^{u-1} e^{-t}\, dt$$
for the incomplete Gamma function. 

\begin{THEB} The following functions are {\rm CM}.

\medskip

{\em (a)} For every $\a\in (0,1],$ the function
$$x\;\mapsto\;\frac{e^x}{\a}\, -\,\Ea(x^\a).$$ 

\smallskip

{\em (b)} For every $\a\in (0,1]$ and $\bt > 1,$ the functions
$$ \frac{x^{1-\bt} e^x\ga(\bt-1,x)}{\a\Ga(\bt -1)}\, -\,\Fab(x)\qquad \mbox{and}\qquad  \Lb\Fab(x)\, -\,\frac{x^{1-\bt} e^x\ga(\bt-1,x)}{\a\Ga(\bt -1)}\cdot$$ 

\smallskip

{\em (c)} For every $\a\in [1,2],$ the function
$$x\;\mapsto\; \Ea(x^\a)\, -\, \frac{e^x}{\a}\cdot$$ 

\smallskip

{\em (d)} For every $\a\in [1,2]$ and $\bt > 1,$ the functions
$$\Fab(x)\, -\, \frac{x^{1-\bt} e^x\ga(\bt-1,x)}{\a\Ga(\bt -1)}\qquad \mbox{and}\qquad   \frac{x^{1-\bt} e^x\ga(\bt-1,x)}{\a\Ga(\bt -1)}\, -\, \Lb\Fab(x).$$ 

\end{THEB} 

\medskip

This result shows that when $\bt \ge 1,$ the CM functions of Theorem A are decomposed, in a non-trivial way, into the sum of two CM functions. It seems however difficult to obtain such a decomposition for $\bt < 1,$ because the underlying Bernstein measures of $\Dab$ or $\cDab$ have then a complicated expression. 

It is interesting to interpret the result of Theorem B in light of the known asymptotic expansions of $\Eab(x^\a)$ and $\ga(\bt-1,x)$ at infinity. Using
$$\int_x^\infty t^{s-1}e^{-t}\, dt\; =\; x^{s-1} e^{-x}\lpa\sum_{k=0}^n \frac{\Ga(s)\,x^{-k}}{\Ga(s-k)}\; +\; o(x^{-(n+1)})\rpa$$
we get the expansion
$$\Fab(x)\, -\, \frac{x^{1-\bt} e^x\ga(\bt-1,x)}{\a\Ga(\bt -1)}\; =\; \sum_{k=1}^n\lpa \frac{x^{-k}}{\a\Ga(\bt-k)}\,-\,\frac{x^{-\a k}}{\Ga(\bt-\a k)}\rpa \; +\; o(x^{-(n+1)(\a\wedge 1)}),$$
whose leading term is negative for $\a < 1$ and positive for $\a > 1,$ in accordance with (b) and (d). Throughout this paper, we will often encounter functions depending on a parameter which become CM when they are positive - see Remark 1 (b) below. Notice also that the  coefficients in the expansion all vanish for $\a =1,$ which matches the easily established - see (\ref{gainc}) below - identity 
$$E_{1,\bt} (x)\; =\; \frac{x^{1-\bt} e^x\ga(\bt-1,x)}{\Ga(\bt -1)}\cdot$$

It is possible to improve Part (a) of Theorem A as well as Parts (a) and (b) of Theorem B, in showing the CM property for the functions taken at $x^{1/\a}.$ For Theorem B the situation is quite different according as $\a < 1/2$ or $\a \in[1/2,1],$ where the functions under consideration have connections with the spectrally positive $(1/\a)-$stable L\'evy process. On the other hand, when $\a\in (1,2]$ the involved functions are linked with the spectrally positive $\a-$stable L\'evy process. All these relationships are explained in details towards the end of the paper. Theorems A and B are proved in Section 2 and 3 respectively, following an approach which is mainly based on Mellin and Stieltjes inversions, and depends in a crucial way on certain Abelian transformations inspired by that of \cite{Mi} and connecting the relevant functions with one another. The argumentation becomes quite intricate for  $0 <\a, \bt < 1$ inside the necktie domain. At the end of each proof, we provide a complete list of the underlying Bernstein measures.

\section{Proof of Theorem A}

\subsection{Proof of (a)} We begin with the only if part. The necessity of $\a\in (0,1)$ comes from the fact that if $\a > 1,$ 
$$\Dab(0+)\; =\; -\frac{1}{\Ga(\a+\bt)}\; < \; 0.$$
Suppose now $\a\in (0,1).$ Then
$$\Dab(x)\; \sim\; \frac{x^{\a -1}}{\Ga(\a+\bt-1)}$$
as $x \to 0+,$ an expression which is negative if $\bt <  1-\a.$ Finally, if $1-\a \le \bt < \a$ then (\ref{CVSR}) entails 
$$\Dab(x)\; =\; x^{\a-1} F_{\a,\a+\bt -1}(x) - \Fab(x)\; \sim\; \frac{x^{-\a}}{\Ga(\bt-\a)}\; < \; 0, \quad x \to +\infty,$$
where the equivalence comes from the aforementioned asymptotic expansion 18.1(22) in \cite{E3}.\\

We now show that if $\a \in (0,1)$ and $\bt \ge \a\vee (1-\a),$ then $\Dab$ is CM. We start with the case $\bt = 1,$ a situation which was already settled in \cite{TS0} - see Remark 2 (b) therein - with the help of Hankel's contour formula for the reciprocal of the Gamma function. We provide here an alternative argument relying on Stieltjes inversion. Using (\ref{CVSR}), we first compute the Laplace transform
$$\int_0^\infty e^{-sx} D_{\a, 1} (x)\, dx \; =\; \frac{1- s^{\a-1}}{s^\a - 1}$$
for every $s > 0.$ The function
$$\Fa(s) \; =\; \frac{1- s^{\a-1}}{s^\a - 1}$$
has an analytic extension on $\CC\slash (-\infty, 0]$ such that $\Fa (z)\to 0$ as $\vert z\vert \to \infty$ and $\Fa (z) = o (\vert z\vert^{-1})$ as $\vert z\vert \to 0,$ uniformly in every sector $\vert \arg z\vert \le \pi -\varepsilon, \,\varepsilon > 0.$ Besides, for every $r > 0$ and $\theta \in (-\pi,\pi)$ one has, setting $\gamma = \a-1/2\in (-1/2,1/2),$
\begin{eqnarray*}
\Re (e^{\i\theta/2} \Fa(re^{\i\theta})) & = & \frac{2(r^{\a-1} +r^\a)\cos (\ga\theta) - 2(1+r^{2\a-1})\cos(\theta/2)}{r^{2\a} - 2r^\a\cos(\theta) +1}\\
& = &  \frac{2(r^{\a-1} +r^\a)(\cos (\ga\theta) -\cos(\theta/2)) + 2(r^\a-1)(1-r^{\a-1})\cos(\theta/2)}{r^{2\a} - 2r^\a\cos(\theta) +1} \; \ge \; 0.
\end{eqnarray*}
From \cite{HW} p. 238, we deduce that $\Fa$ is the Stieltjes transform of some positive measure $\mu_\a (dt)$ on $(0,+\infty),$ viz. it writes
$$\Fa (s)\; =\; \int_0^\infty \frac{\mu_\a(dt)}{t+s}\; =\; \int_0^\infty e^{-sx} \lpa \int_0^\infty e^{-xt} \mu_\a(dt)\rpa dx.$$
Moreover, the Perron-Stieltjes inversion formula - see e.g. Theorem VIII.7.a p. 339 in \cite{W} - entails 
\begin{eqnarray*}
\mu_\a(b)\, -\, \mu_\a(a) & = & \lim_{\eta \to 0+} \frac{1}{2\pi\i }\int_a^b (\Fa(-t-\i\eta) - \Fa(-t+\i\eta))\, dt\\
& = & \int_a^b \frac{\sin(\pi\a) t^{\a -1} (1+t)}{\pi(t^{2\a} - 2\cos(\pi\a)t^\a +1)}\, dt
\end{eqnarray*}
for every $0< a <b.$ By uniqueness of the Laplace transform, this yields
\begin{equation}
\label{b1}
D_{\a,1}(x)\; =\; \int_0^\infty e^{-xt}  \frac{\sin(\pi\a) t^{\a -1} (1+t)}{\pi(t^{2\a} - 2\cos(\pi\a)t^\a +1)}\, dt, \qquad x >0,
\end{equation}
and shows that $D_{\a,1}$ is CM. To handle the case $\bt > 1$ we appeal to the formula
\begin{equation}
\label{Frac1}
D_{\a,\bt}(x)\; =\; \frac{1}{\Ga(\bt -1)}\int_0^1 (1-t)^{\bt -2} D_{\a, 1} (xt) \, dt,
\end{equation}
which can be checked from (\ref{CVSR}). Setting
$$\fa(t)\; =\; \frac{\sin(\pi\a) t^{\a -1} (1+t)}{\pi(t^{2\a} - 2\cos(\pi\a)t^\a +1)}$$
we deduce from (\ref{Frac1}) and Fubini's theorem
$$D_{\a,\bt}(x)\; =\; \int_0^\infty e^{-xt} \lpa \frac{1}{\Ga(\bt -1)}\int_0^1 (1-u)^{\bt -2} \fa (\frac{t}{u}) \, \frac{du}{u}\rpa dt, \qquad x >0.$$
This shows that $\Dab$ is CM for every $\bt > 1.$ \\

We now proceed to the case $\a\vee(1-\a) \le \bt < 1$ which is more delicate, except in the case $\a = 1/2$ where  (\ref{CVSR}) entails
$$D_{1/2, \bt} (x)\; =\; \frac{x^{-1/2}}{\Ga(\bt-1/2)},$$
a CM function if $\bt \ge 1/2$ (and the zero function if $\bt = 1/2$). We must divide the proof according as $\a > 1/2$ or $\a < 1/2.$ 

\subsubsection{The case $\a > 1/2$ and $\bt\ge \a$} A computation based on (\ref{CVSR}) reveals that
\begin{equation}
\label{12+}
\Dab(x) \; = \; \frac{1}{\Ga(\bt -\a)}\int_0^1 (1-t)^{\bt -\a -1} t^{\a -1} D_{\a, \a} (xt) \, dt
\end{equation}
for every $\bt > \a,$ so that it is enough to show that $D_{\a,\a}$ is CM. The proof of this fact hinges upon a certain multiplicative factorization of the $\sigma-$finite measure $\mu_\a$ on $(0,+\infty)$ with density $\fa.$ The Mellin transform\footnote{Throughout this paper we integrate along $t^s$ instead of $t^{s-1}$ to define Mellin transforms, because this leads to shorter formul\ae.}
$$\Ma(s)\; =\; \int_0^\infty t^s \fa(t)\, dt$$
is well-defined for $s\in(-\a, \a -1)$ and a computation similar to Proposition 4 in \cite{TS1} yields the closed formula
$$\Ma(s) \; = \; - \frac{\sin(\pi/\a)\sin(\pi s)}{\a \sin (\pi s /\a)\sin (\pi(s+1)/\a)}\cdot$$
Introduce the $\B_{\a, 1-\a}$ random variable with density
$$\frac{x^{\a-1}(1-x)^{-\a}}{\Ga(\a)\Ga(1-\a)}\Un_{(0,1)}(x)$$
and with Mellin transform
$$\EE[\B_{\a, 1-\a}^s]\; =\; \frac{\Ga(\a +s)}{\Ga(1+s)\Ga(\a)}, \qquad s > -\a.$$
The complement and concatenation formul\ae\, for the Gamma function entail
$$\frac{\Ga(1+s)\Ga(\a)\Ma(s)}{\Ga(\a +s)}\; =\; \frac{\Ga(1-s/\a)}{\Ga(1-s)}\,\times\,  \frac{\Ga(1+s/\a)\Ga(\a)}{\Ga(\a+s)}\,\times\,  \frac{-\sin(\pi/\a)}{\sin(\pi(s+1)/\a)}$$
for $s\in(-\a, \a -1).$ We will now show that the three factors on the right hand side are Mellin transforms of positive $\sigma-$finite measures on $(0,+\infty).$ First, recall e.g. from Theorem 2.6.3 in \cite{Z} that
$$\frac{\Ga(1-s/\a)}{\Ga(1-s)}\; =\; \EE[\Za^s], \qquad s < \a,$$
where $\Za$ is the standard positive $\a-$stable random variable which is defined through the Laplace transform
$$\EE[e^{-\lbd \Za}]\; =\; e^{-\lbd^\a}, \qquad \lbd \ge 0.$$
The following Lemma, which might be well-known although we could locate it in the literature, shows the property for the second factor.

\begin{LEM} For every $\a\in (0,1)$ there exists a positive random variable $\Xa$ such that
$$\frac{\Ga(1+s/\a)\Ga(\a)}{\Ga(\a+s)}\; =\; \EE[\Xa^s], \qquad s > -\a.$$
\end{LEM}

\proof A simple application of Helly's selection theorem, and the fact that the quotient on the left-hand side equals 1 at $s =0,$ show that it is enough to consider the case when $\a$ is rational. Set $\a = p/q$ with $q > p\ge 1$. A repeated use of the Legendre-Gauss multiplication formula for the Gamma function entails
\begin{eqnarray*}
\frac{\Ga(1+\frac{qs}{p})\Ga(\frac{p}{q})}{\Ga(\frac{p}{q}+s)} & = & \lpa \frac{q^{\frac{q}{p}}}{p^{\frac{p}{q}}}\rpa^s\prod_{i=2}^{p} \frac{\Ga(\frac{i}{q}+\frac{s}{p})\Ga(\frac{i-1}{p} + \frac{1}{q})}{\Ga(\frac{1}{q}+\frac{i-1}{p} +\frac{s}{p})\Ga(\frac{i}{q})} \;\,\times\prod_{j=p+1}^{q} \frac{\Ga(\frac{j}{q}+\frac{s}{p})}{\Ga(\frac{j}{q})}\;\; =\;\; \EE[{\bf X}^s_{p,q}]
\end{eqnarray*}
where $\B_{a,b}$ and $\G_c$ stand for the Beta and Gamma random variables with respective parameters $a,b,c >0,$ and ${\bf X}_{p,q}$ is the independent product
$$\frac{q^{\frac{q}{p}}}{p^{\frac{p}{q}}}\;\,\times \lpa\prod_{i=2}^{p} \B_{\frac{i}{q}, (i-1)(\frac{1}{p} -\frac{1}{q})} \;\,\times\prod_{j=p+1}^{q} \G_{\frac{j}{q}}\rpa^{\frac{1}{p}}\cdot$$
This completes the proof. 

\endproof
To handle the third factor we rewrite, for every $s\in(-\a,\a-1),$
\begin{eqnarray*}
\frac{-\sin(\pi/\a)}{\sin(\pi(s+1)/\a)} & = &  \frac{\Ga((1+s)/\a)}{\Ga(1/\a)}\; \times\;\frac{\Ga(1- (1+s)/\a)}{-\Ga(1- 1/\a)}\\
& = & \EE[\G_{1/\a}^{s/\a}]\; \times\;\int_0^\infty\lpa \frac{(1-\a)t^{-\a}e^{-t^{-\a}}}{\Ga(2-1/\a)}\rpa t^s \, dt.
\end{eqnarray*}
Setting $\gaa$ for the integrated function on the right-hand side and putting everything together shows finally that
$$\Ma(s)\; =\; \EE[\B_{\a, 1-\a}^s]\,\times\,\EE[\Ya^s]\,\times\,\M_{\gaa} (s), \qquad s\in(-\a, \a-1),$$
where $\Ya = \Za\times \Xa \times \G_{\frac{1}{\a}}^{\frac{1}{\a}}$ is meant as an independent product and we denote by
$$\M_f(s)\; =\; \int_0^\infty f(t)\, t^s\, dt$$ 
the Mellin transform of a positive measurable function on $(0,+\infty).$ Set now $\odot$ for the multiplicative convolution of two positive measurable functions on $(0,+\infty)$ :
$$f\odot g(x)\; =\; \int_0^\infty f(t)\, g(\frac{x}{t})\, \frac{dt}{t}, $$
and recall from Fubini's theorem that $\M_{f\odot g}\, =\, \M_f\, \times\, \M_g$ (with possible infinite values). Setting $f_\X$ for the density of an absolutely continuous random variable $\X$ and $\ha = \gaa \odot f_{\Ya},$ we get

$$\Ma(s)\; =\; \M_{f^{}_{\B_{\a, 1-\a}}}\!(s)\;\times\;\M_{\ha}(s)\; <\; +\infty$$
for every $s\in (-\a,\a-1).$ Inverting these Mellin transforms entails the crucial factorization
\begin{equation}
\label{Mell}
\fa\; =\; f^{}_{\B_{\a, 1-\a}}\!\odot\;\ha, 
\end{equation}
which seems difficult to obtain from a direct computation. We can now finish the proof of the case $\a > 1/2.$ On the one hand, it follows from (\ref{12+}) with $\bt = 1$ that
$$D_{\a, 1}(x)\; =\; \Ga(\a)\int_0^1 f^{}_{\B_{\a, 1-\a}}(t)\, D_{\a,\a}(xt)\, dt.$$
On the other hand, we deduce from (\ref{b1}), (\ref{Mell}) and Fubini's theorem
$$D_{\a, 1}(x)\; =\; \int_0^1 f^{}_{\B_{\a, 1-\a}}(t)\lpa\int_0^\infty e^{-xtu} \ha(u)\, du\rpa dt.$$
Setting
$$\Ha(x) \; =\; \frac{1}{\Ga(\a)}\int_0^\infty e^{-xu} \ha(u)\, du, \qquad x >0,$$
which is a CM function, and comparing the two identities, we obtain
$$\int_0^1 f^{}_{\B_{\a, 1-\a}}(t)\, D_{\a,\a}(xt)\, dt\; =\; \int_0^1 f^{}_{\B_{\a, 1-\a}}(t)\, \Ha(xt)\, dt, \qquad x > 0.$$
This is an identity between certain Abelian transforms for which we did not find any direct inversion formula in the literature. After a change of variable, this identity changes into 
$$\ka\,\odot\, D_{\a,\a}\; =\; \ka\, \odot\, \Ha$$
with $\ka(x) = (x-1)^{-\a}\Un_{\{x > 1\}}.$ The Mellin transform of the left-hand side is well-defined on $(-\a,\a -1)$ because
$$D_{\a,\a}(x)\, \sim\, \frac{x^{\a-1}}{\Ga(2\a-1)}\quad\mbox{as $x\to 0+$}\qquad\mbox{and}\qquad D_{\a,\a}(x)\, \sim\, -\frac{x^{-1}}{\Ga(\a-1)}\quad\mbox{as $x\to +\infty.$}$$
Hence, taking the Mellin transform on both sides and factorizing by that of $\ka$ entails finally, after Mellin inversion, that $D_{\a,\a} = \Ha$ and the proof is complete.

\qed

\subsubsection{The case $\a < 1/2$ and $\bt\ge 1-\a$} Another direct computation based on (\ref{CVSR}) yields
\begin{equation}
\label{12-}
\Dab(x) \; = \; \frac{x^{\a-1}}{\Ga(\a+\bt -1)} \; + \; \frac{1}{\Ga(\a+\bt -1)}\int_0^1 (1-t)^{\a+\bt -2} t^{-\a} D_{\a, 1-\a} (xt) \, dt
\end{equation}
for every $\bt > 1-\a,$ so that it is enough to show that $D_{\a,1-\a}$ is CM. As before, this fact will follow from a certain multiplicative factorization of the measure $\mu_\a.$ However the Mellin transform of $\mu_\a$ is here everywhere infinite,  because $\a<1/2$. Set $\na = [\a^{-1}] \ge 2,$ which is the unique integer such that
$$\frac{1}{\na+1}\; < \; \a\; \le\;\frac{1}{\na}\cdot$$  
We will show the stronger result that the function
\begin{equation}
\label{till}
{\widetilde D_{\a, 1-\a}} (x)\; =\; D_{\a, 1-\a} (x) \; -\;  \sum^{\na-1}_{k=2}\, \frac{x^{\a k-1}}{\Ga (\a (k-1))}
\end{equation}
(with an empty sum if $\na =2$) is CM. This result is indeed stronger because the sum on the right-hand side is clearly CM. Observe first from (\ref{CVSR}) that ${\widetilde D_{\a, 1-\a}} = 0$ if $\a$ is the reciprocal of an integer viz. $\a = 1/\na.$ From now on we will hence suppose $\a < 1/\na.$  Decompose
\begin{eqnarray*}
D_{\a, 1} (x) &  = &\sum^{\na-1}_{k=1}\, \frac{x^{\a k-1}}{\Ga (\a k)}\; +\;\sum_{k\ge \na}\, \frac{x^{\a k-1}}{\Ga (\a k)}\; -\;\sum_{k\ge 0}\, \frac{x^{\a k}}{\Ga (1+\a k)} \\
& = & \sum^{\na-1}_{k=1}\, \frac{\sin(\pi\a k)}{\pi}\int_0^\infty e^{-xt} t^{-\a k} dt \; +\;\sum_{k\ge \na}\, \frac{x^{\a k-1}}{\Ga (\a k)}\; -\;\sum_{k\ge 0}\, \frac{x^{\a k}}{\Ga (1+\a k)} \\
& = & \frac{\sin(\pi\a)}{\pi}\int_0^\infty\!\! e^{-xt} \lpa \sum^{\na-1}_{k=1} U_{k-1}(\cos\pi\a)\, t^{-\a k}\rpa dt \, +\sum_{k\ge \na} \frac{x^{\a k-1}}{\Ga (\a k)}\, -\,\sum_{k\ge 0}\, \frac{x^{\a k}}{\Ga (1+\a k)} 
\end{eqnarray*}
where 
$$U_n(\cos\theta)\; =\; \frac{\sin(n+1)\theta}{\sin\theta}$$
stands for the $n-$th Chebyshev polynomial of the second kind. Using the notation $U_n^\a = U_n(\cos\pi\a)$ for every $n\ge 0$ and 
$${\widetilde D_{\a, 1}} (x)\; =\; \sum_{k\ge \na}\, \frac{x^{\a k-1}}{\Ga (\a k)}\; -\;\sum_{k\ge 0}\, \frac{x^{\a k}}{\Ga (1+\a k)},$$ 
we get
$${\widetilde D_{\a, 1}} (x)\; =\; \int_0^\infty  e^{-xt} \tfa (t)\, dt$$
where, simplifying with the help of the recurrence relations $U_{n+2} + U_n = U_1U_{n+1},$ we compute
\begin{eqnarray*}
 \tfa (t) & = & \frac{\sin(\pi\a) t^{\a -1} (1+t)}{\pi(t^{2\a} - U_1^\a t^\a +1)} \; -\;  \frac{\sin(\pi\a)}{\pi}\, \sum^{\na-1}_{k=1} U_{k-1}^\a t^{-\a k}\\
& = & \frac{\sin(\pi\a) (t^{\a -1}  - U_{\na -2}^\a t^{-\a(\na-1)}+ U_{\na -1}^\a t^{-\a(\na-2)})}{\pi(t^{2\a} - U_1^\a t^\a +1)}\cdot
\end{eqnarray*}
From (\ref{12-}) with $\bt =1$ and the fact that
$$\sum^{\na-1}_{k=2}\, \frac{x^{\a k-1}}{\Ga (\a k)}\; =\; \frac{1}{\Ga(\a)}\int_0^1 (1-t)^{\a-1} t^{-\a} \lpa\sum^{\na-1}_{k=2}\, \frac{(xt)^{\a k-1}}{\Ga (\a (k-1))}\rpa \, dt$$
we also have
\begin{equation}
\label{DD}
\int_0^\infty  e^{-xt} \tfa (t)\, dt \; = \;  \frac{1}{\Ga(\a)}\int_0^1 (1-t)^{\a-1} t^{-\a} {\widetilde D_{\a, 1-\a}} (xt) \, dt.
\end{equation}
Similarly as in the case $\a > 1/2,$ our next step is now to show that the density of the random variable $\B_{1-\a,\a}$ is a multiplicative factor of $\tfa.$ The Mellin transform
$$\cMa (s)\; =\; \int_0^\infty t^s \tfa (t)\, dt$$
is finite for every $s\in (-\a, \na\a - 1)$ which is a non-empty interval. Using
$$\frac{\sin\pi\a}{\pi} \int_0^\infty \frac{t^{\a -1+s} dt}{t^{2\a} - U_1^\a t^\a +1} \; =\; \frac{\sin\pi\a}{\pi\a} \int_0^\infty \frac{u^{\frac{s}{\a}}du}{u^2 + 2 u\cos\pi(1-\a) +1} \; =\; \frac{\sin(\pi s(1-\a)/\a)}{\a \sin(\pi s/\a)}$$
for every $s\in (-\a, \a),$ where the second computation comes from a standard application of the residue theorem, we deduce 
\begin{eqnarray*}
\cMa (s) & = & \frac{\sin(\pi s(1-\a)/\a)}{\a \sin(\pi s/\a)}\; -\; \frac{U_{\na -2}^\a \sin(\pi (1-\a)(s +1-\na\a)/\a)}{\a \sin(\pi (s +1-\na\a)/\a)}\\
& & \qquad +\;\, \frac{U_{\na -1}^\a\sin(\pi (1-\a)(s +1+\a -\na\a)/\a)/\a)}{\a \sin(\pi (s +1+\a-\na\a)/\a)/\a)}\\
& = & - \frac{\sin(\pi/\a)\sin(\pi s)}{\a \sin (\pi s /\a)\sin (\pi(s+1)/\a)}
\end{eqnarray*}
 for every $s\in (-\a, \na\a - 1),$ where the second equality follows after some trigonometry. This entails, with the above notation,
\begin{eqnarray*}
\frac{\Ga(1+s)\Ga(1-\a)\cMa(s)}{\Ga(1-\a +s)} & = & \frac{\Ga(1-s/\a)}{\Ga(1-s)}\;\times\;\frac{\Ga(1+s/\a)\Ga(\a)}{\Ga(\a+s)} \\
& & \qquad \times\;  \frac{\Ga(\a+s)\Ga(1-\a)}{\Ga(1-\a+s)\Ga(\a)}\;\times\;  \frac{-\sin(\pi/\a)}{\sin(\pi(s+1)/\a)}\\
& = & \M_{f^{}_{\Za}}\!(s)\;\times\;\M_{f^{}_{\Xa}}\!(s)\;\times\;\M_{f^{}_{\B_{\a, 1-2\a}}}\!(s)\;\times\;\frac{-\sin(\pi/\a)}{\sin(\pi(s+1)/\a)}\cdot
\end{eqnarray*}
Using the concatenation formula for the Gamma function and setting $\U$ for the uniform random variable on $(0,1),$ we finally decompose
\begin{eqnarray*}
\frac{-\sin(\pi/\a)}{\sin(\pi(s+1)/\a)} & = &  \frac{\Ga((1+s)/\a)}{\Ga(1/\a)}\; \times\;\frac{\Ga(1- (1+s)/\a)}{-\Ga(1- 1/\a)}\\
& = & \EE[\G_{1/\a}^{s/\a}]\; \times\;\prod_{k=1}^{\na -1} \lpa\frac{1-\a k}{1-\a k + s}\rpa\; \times\;\frac{\Ga(\na - 1/\a -s/\a)(1/\a -\na)}{\Ga(\na +1 - 1/\a)}\\
& = & \EE[\G_{1/\a}^{s/\a}]\; \times\;\prod_{k=1}^{\na -1} \EE[\U^{s/(1-\a k)}] \; \times\;\int_0^\infty\lpa \frac{(1-\na\a)t^{-\na\a}e^{-t^{-\a}}}{\Ga(\na +1-1/\a)}\rpa t^s \, dt
\end{eqnarray*}
for every $s\in (-\a, \na\a -1).$ Setting $\tga$ for the integrated function on the right-hand side and
$$\tha \; =\; \tga\;\odot\; f^{}_{\B_{\a, 1-2\a}}\;\odot\;f^{}_{\Za}\;\odot\;f^{}_{\Xa}\;\odot\;f^{}_{\G^{1/a}_{1/\a}}\;\odot\;\lpa \bigodot_{k=1}^{\na -1} f^{}_{\U^{1/(1-\a k)}}\rpa,$$
an inversion of the Mellin transform yields the factorization
$$\tfa\; =\; f^{}_{\B_{1-\a, \a}}\;\odot\; \tha.$$
Comparing with (\ref{DD}), we deduce
$$\int_0^1  f^{}_{\B_{1-\a, \a}} {\widetilde D_{\a, 1-\a}} (xt) \, dt\; =\; \int_0^1  f^{}_{\B_{1-\a, \a}} \tHa (xt) \, dt, \qquad x >0,$$
where $\tHa$ is the CM function
$$\tHa(x)\; =\; \frac{1}{\Ga(1-\a)}\int_0^\infty e^{-xt}\, \tha (t)\, dt, \qquad x > 0.$$
The latter identity transforms into 
$$\tka\,\odot\, {\widetilde D_{\a,1-\a}}\; =\; \tka\, \odot\, \tHa$$
with $\tka(x) = (x-1)^{\a-1}\Un_{\{x > 1\}}.$ The Mellin transform of the left-hand side is well-defined on the non-empty interval $(-\na\a,\a -1)$ because
$${\widetilde D_{\a,1-\a}}(x)\, \sim\, \frac{x^{\na\a-1}}{\Ga(\na\a-\a)}\;\;\mbox{as $x\to 0+$}\qquad\mbox{and}\qquad {\widetilde D_{\a,1-\a}} (x)\, \sim\, \frac{x^{-\a}}{\Ga(1-2\a)}\;\;\mbox{as $x\to +\infty.$}$$
Similarly as above, we obtain the identification ${\widetilde D_{\a,1-\a}} = \tHa$ and the proof is complete.

\qed

\subsection{Proof of (b)} We begin with the only if part, which is analogous to the above. The necessity of $\a\le 2$ comes from the fact that if $\a  >2,$ then $\cDab(0+)' =0$ so that $\cDab$ is not CM. Suppose now $\a\in (0,1).$ Then
$$\cDab(x)\; \sim\; \frac{-x^{\a -1}}{\Ga(\a+\bt-1)}\quad\mbox{as $x \to 0+$}\qquad\mbox{and}\qquad \cDab(x)\; \sim\; \frac{-x^{-\a}}{\Ga(\bt-\a)}\quad\mbox{as $x \to +\infty,$}$$ 
and at least one of these expressions is negative if $\bt >  \a\wedge(1-\a).$ Finally, if $\a\in(1,2]$ then again the expansion 18.1.(22) in \cite{E3} entails 
$$\cDab(x)\; \sim\; \frac{x^{-1}}{\Ga(\bt-1)}\quad \mbox{as $x \to +\infty,$}$$
which is negative if $\bt < 1.$ \\

We next show that if $\a \in (1,2]$ and $\bt \ge 1,$ then $\cDab$ is CM. The case $\bt = 1$ is stated as Theorem 1 in \cite{TS0} and can also be handled with exactly the same Stieltjes inversion argument as in the proof of (a). For $\a < 2$ this reads
\begin{equation}
\label{b2}
\overline{D_{\a,1}}(x)\; =\; \int_0^\infty e^{-xt}  \frac{-\sin(\pi\a) t^{\a -1} (1+t)}{\pi(t^{2\a} - 2\cos(\pi\a)t^\a +1)}\, dt, \qquad x >0,
\end{equation}
whereas for $\a = 2$ we simply have $\overline{D_{2,1}}(x) = e^{-x},$ the prototype of a CM function. The case $\bt > 1$ follows from the formula
\begin{equation}
\label{b3}
\cDab(x)\; =\; \frac{1}{\Ga(\bt -1)}\int_0^1 (1-t)^{\bt -2} \overline{D_{\a,1}}(xt)\, dt
\end{equation}
which is again a direct computation relying on (\ref{CVSR}).\\

We finally show that if $\a\in (0,1)$ and $\bt \le \a\wedge (1-\a),$ then $\cDab$ is CM. This is the delicate part. After some further computations relying on (\ref{CVSR}), we get 
\begin{eqnarray*}
\cDab(x) \; = \; x^{1-\bt} \frac{d}{dx}(x^\bt \overline{D_{\a, \bt +1}}(x)) & = & -(xD_{\a, \bt +1}'(x) +\bt D_{\a, \bt +1}(x))\\
& = & \frac{1}{\Ga(\bt)}\int_0^1 (1-t)^{\bt -1} \Hab(xt)\, dt
\end{eqnarray*}
with
\begin{eqnarray*}
\Hab (x) & = & -(\bt D_{\a, 1}(x) + x D_{\a, 1}'(x))\; = \; \int_0^\infty e^{-xt} (t\fa'(t)+(1-\bt)\fa(t))\, dt
\end{eqnarray*}
where the last equality comes after an integration by parts in (\ref{b1}). We are hence reduced to show that $\Hab$ is CM, in other words that the function 
$$\hab (t) \; =\; (1-\bt)\fa(t) +t\fa'(t)$$
is non-negative. We must divide the proof according as $\a \le 1/2$ or $\a > 1/2.$\\

(i) The case $\a \le 1/2$ and $\bt\le \a$. It is enough to show the non-negativity of $h_{\a, \a}$ which, after some computations, amounts to that of the function
$$t\mapsto(1-2\a)t^{3\a} - 2(1-\a) \cos(\pi\a) t^{2\a} + t^\a - 2\a t^{3\a -1} + 2\a \cos (\pi\a) t^{2\a-1}.$$

(ii) The case $\a > 1/2$ and $\bt\le 1-\a$. Here we need to show the non-negativity of $h_{\a, 1-\a},$ which is equivalent to that of the function
$$t\;\mapsto\; -2\a\cos(\pi\a) t^{2\a} - t^{3\a-1} + 2\a t^\a + 2(1-\a)\cos(\pi\a) t^{2\a -1} + (2\a -1) t^{\a-1}.$$

The non-negativity of these two polynomial functions, which are zero for $\a =1/2,$ can be observed heuristically with the help of some plotting software. A strict proof 
is obtained in using the sequence of signs $+-+-+$ for the coefficients, and the following equivalence for all $a,b,c,\rho >0:$ 
$$at^{1+\rho} -bt + c\, \ge\, 0\;\,\mbox{for all}\; t > 0\quad \Leftrightarrow\quad a\, \ge \,\lpa \frac{\rho}{c}\rpa^\rho \lpa \frac{b}{1+\rho}\rpa^{1+\rho}.$$
Though interesting, the details of this strict proof are lenghty and will be not included here. They have been typesetted and are available upon request. 

\qed

\begin{REMS} \label{Sab} {\em (a) Proceeding along the same way as in \cite{TS0}, we obtain the following integral representation which is valid for any $0 <\a,\bt < 2.$
$$\cDab (x) \; =\; \frac{x^{1-\bt}}{2\pi \i}\int_H \lpa \frac{t^{\a -1} -1}{t^\a -1}\rpa t^{1-\bt} e^{xt} \,dt,$$
where $H$ is a standard Hankel path encircling the origin. This transforms into
$$\frac{x^{1-\bt}}{\pi}\int_0^\infty e^{-xt}\lpa \frac{\sin(\pi\bt) t^{1-\bt}(t^{2\a -1} -1) + \sin(\pi(\a+\bt)) t^{\a +1-\bt} + \sin(\pi(\a-\bt)) t^{\a-\bt}}{t^{2\a} - 2\cos(\pi\a)t^\a +1}\rpa   dt$$
and simplifies to (\ref{b1}) resp. (\ref{b2}) for $\bt = 1.$ Using this representation it is an easy exercise, which is left to the reader, to prove the weaker result that for any $\a\in (0,1)$ the function 
$$x^{\bt -1}\Dab(x)\qquad\mbox{resp.}\qquad x^{\bt -1}\cDab(x)$$ 
is CM when $\bt\in [\a\vee (1-\a), 1)$ resp. $\bt\in (0,\a\wedge (1-\a)].$ \\

(b) It is clear from the above proofs that for any $\a\in (0,2]$ the following equivalences hold
$$\Dab \;\,\mbox{is CM}\;\Leftrightarrow\; \Dab\; \ge\; 0\qquad\mbox{and}\qquad\cDab \;\,\mbox{is CM}\;\Leftrightarrow\; \cDab\; \ge\; 0.$$
When $\a > 4,$ the asymptotic expansion 18.1 (22) p.210 in \cite{E3} shows that $\cDab$ has leading term
$$4 \cos (2\pi (\a-1)/\a) e^{\cos (2\pi/\a) x} x^{1-\bt} \cos (x \sin (2\pi/\a) + 2\pi(3-\a-\bt)/\a)$$
at infinity, and hence oscillates. It can also be shown directly that 
$$\overline{D_{4,1}}(x) \; =\; \frac{1}{2}\, (e^{-x} + \cos (x) + \sin(x))$$
and also oscillates. It is however unclear to the author whether the function $\cDab$ is everywhere non-negative for $\a\in (2,4], \bt >0,$ or not. Indeed, the leading term at infinity is then a positive monomial. See \cite{WZ} for more complete results on the asymptotic expansions of Mittag-Leffler functions. }

\end{REMS}

\subsection{List of the Bernstein measures and an improvement of (a)}

In this paragraph we recapitulate the explicit Bernstein measures $\mab$ associated with the CM functions $\vert\Dab\vert.$ Recall that these are the positive $\sigma-$finite measures such that
$$\vert\Dab\vert (x)\; =\; \int_0^\infty e^{-xt} \mab (dt)\, dt, \qquad x > 0.$$
For the sake of completeness we also rephrase the Bernstein measures associated with the main result of \cite{Sch}. Finally, we improve part (a) and show that the function 
$$x\;\mapsto\;\Dab (x^{1/\a})$$
is also CM for every $\a\in (0,1)$ and $\bt\ge \a\vee (1-\a).$

\subsubsection{The case $\a\in (1,2]$} In this case one has $\cDab (0+) =1/\Ga(\bt)$ and the underlying Bernstein measure is hence finite. For $\bt =1$ it follows from (\ref{b2}) that
$$\overline{D_{\a, 1}}(x)\; =\; \EE[e^{-x \Ta}]$$
where ${\bf T}_2 = \Un$ and $\Ta$ is for $\a\in (1,2)$ the positive random variable with density
$$\frac{-\sin(\pi\a) t^{\a -1} (1+t)}{\pi(t^{2\a} - 2\cos(\pi\a)t^\a +1)}\cdot$$
This latter random variable is connected to first-passage time of spectrally positive $\a-$stable L\'evy processes - see Theorem 3 in \cite{TS0}, a connection which will be discussed in further details in the next section. For $\bt > 1,$ a direct consequence of (\ref{b3}) and Fubini's theorem is that
\begin{equation}
\label{Sabb}
\cDab(x)\; =\; \frac{1}{\Ga(\bt)}\; \EE[e^{-x (\B_{1,\bt -1}\times\Ta})],
\end{equation}
where here and throughout the product is assumed to be independent. Notice that for $\a = 2$ one obtains
$$\overline{D_{2, \bt}}(x)\; =\; \frac{1}{\Ga(\bt)}\; \EE[e^{-x \B_{1,\bt -1}}],$$
in accordance with (2.4) in \cite{Sch} because 
$$\overline{D_{2, \bt}}(x)\; =\; \sum_{n\ge 0} \frac{(-x)^n}{\Ga(\bt +n)}\; =\; E_{\bt, 1}(-x).$$
\begin{REMS}
\label{ijj}
{\em Using the representation of Remark \ref{Sab} (a) and some computations, we also obtain the explicit formula
\begin{eqnarray*}
\overline{D_{3/2,3/2}}(x) & = & \frac{1}{\pi\sqrt{x}}\int_0^\infty e^{-xt}\lpa \frac{1-t^2}{\sqrt{t}(t^3 +1)}\rpa   dt\\
& = & \frac{1}{\pi\sqrt{\pi}}\int_0^\infty e^{-xt}\lpa \int_0^t \frac{(1-s) ds}{\sqrt{s(t-s)}(s^2 -s+1)}\rpa   dt\\
& = & \frac{2}{\sqrt{3\pi}}\int_0^\infty e^{-xt}\lpa \frac{\Re (\sqrt{{\rm j} -t})}{\sqrt{t^2 - t +1}}\rpa   dt
\end{eqnarray*}
with the notation ${\rm j} = e^{2\i\pi/3}.$ Observe that the function between brackets on the third line is indeed positive. This formula can also be obtained directly from (\ref{Sabb}). Except in this particular case $\a=\bt = 3/2,$ it does not seem that anything more explicit can be obtained from the integral representation of Remark \ref{Sab} (a).}
\end{REMS}

\subsubsection{The case $\a\in (0,1)$} In this case the underlying Bernstein measure is infinite since $\vert \Dab\vert (0+) = +\infty.$ It also follows from the above proofs that this measure has a density, which we denote by $\dab.$ The latter will be expressed with the notations introduced during the above proofs, to which we refer without further repetition. We distinguish several subcases.\\

\begin{itemize}

\item $\bt \ge 1.$ It follows from (\ref{b1}) and (\ref{Frac1}) that $d_{\a, 1} = \fa$ and
that for every $\bt > 1,$
\begin{equation}
\label{dacis}
\dab\; =\; \frac{1}{\Ga(\bt)} \; f_{\B_{1,\bt-1}^{}}\odot\, \fa.
\end{equation}

\item $\a > 1/2$ and $\bt\in [\a, 1).$ One has
\begin{equation}
\label{daa}
d_{\a, \a} \;=\;\frac{1}{\Ga(\a)} \;\ha\; =\; \frac{1}{\Ga(\a)} \; f_{\Za}\,\odot\,f_{\Xa}\,\odot\,f_{\G_{1/\a}^{1/\a}}\,\odot\, \gaa
\end{equation}
and, for every $\bt > \a,$
$$\dab\; =\; \frac{1}{\Ga(\bt)} \; f_{\B_{\a,\bt-\a}^{}}\odot\, \ha.$$

\item $\a \le 1/2$ and $\bt\in [1-\a, 1).$ It follows from (\ref{till}) that
\begin{equation}
\label{da1a}
d_{\a, 1-\a}(x) \;=\;\frac{1}{\Ga(1-\a)} \;\tha(x)\; +\;\sum^{\na-1}_{k=2}\, \frac{x^{-\a k}}{\Ga(1-\a k)\Ga (\a (k-1))}\cdot
\end{equation}
Recall that $\tha$ is zero if $\a = 1/\na$ viz. $\a$ is the reciprocal of an integer, and that the sum is empty if $\na = 2$ viz. $\a =1/2.$ In particular one has $d_{1/2,1/2} =0,$ in accordance with $D_{1/2,1/2} = 0.$ For $\bt > 1-\a,$ from (\ref{12-}) and (\ref{till}) we deduce
$$\;\qquad\Dab(x) \; = \; \frac{1}{\Ga(\a+\bt -1)}\int_0^1 (1-t)^{\a+\bt -2} t^{-\a} {\widetilde D}_{\a, 1-\a} (xt) \, dt\; + \;\sum^{\na-1}_{k=1} \frac{x^{\a k-1}}{\Ga(\a k+\bt -1)},$$
so that 
$$\dab(x) \;=\;\frac{1}{\Ga(\bt)} \;f_{\B_{1-\a,\a+\bt-1}^{}}\!\odot\,\tha(x)\; +\;\sum^{\na-1}_{k=1}\, \frac{x^{-\a k}}{\Ga(1-\a k)\Ga(\a k+\bt -1)}\cdot$$

\item $\a \le 1/2$ and $\bt \le \a\wedge (1-\a).$ The Bernstein density is
$$\dab\; =\;\frac{1}{\Ga(\bt +1)} \;f_{\B_{1,\bt}^{}}\!\odot\,\hab$$
where
$$\hab(t)\; = \; (1-\bt)\fa(t) + t\fa'(t)\; = \;  \frac{\sin(\pi\a)\thab (t)}{\pi(t^{2\a} - 2\cos(\pi\a)t^\a +1)^2}$$
and we have set

\begin{eqnarray*}
\thab (t) & = & (1-\a-\bt)t^{3\a}\;-\;2(1-\bt)\cos(\pi\a) t^{2\a} \;+ \;(1-\a-\bt) t^\a \\
& & \qquad\qquad -\; (\a +\bt) t^{3\a-1} \; +\; 2\bt\cos(\pi\a) t^{2\a -1}\; +\; (\a -\bt) t^{\a-1}.
\end{eqnarray*}
Observe that this latter function $\thab$ takes negative values when $\bt > \a\wedge (1-\a).$

\end{itemize}

\subsubsection{Bernstein densities associated with $\Eab(-x)$} In this paragraph we express the positive random variables underlying the CM functions $\Eab(-x)$ for $\a\in (0,1), \bt \ge \a.$ This is basically a reformulation of the main results of \cite{Sch, Mi}, which we however believe to be worth mentioning. In the case $\a = 1,$ we saw above during the analysis of $D_{2, \bt}$ that 
$$E_{1,\bt}(-x)\; =\; \frac{1}{\Ga(\bt)}\;\EE[e^{-x\B_{1, \bt -1}}].$$ 
The classical case $\bt = 1$ had been settled in \cite{P}, whose main result reads
$$\Ea(-x)\; =\; \EE[e^{-x\MMa}]$$
where $\MMa \elaw \Za^{-\a}$ is the so-called Mittag-Leffler random variable. It follows from Lemma 1 in \cite{Mi} that
$$E_{\a,\a}(-x)\; =\; \a \EE[\MMa e^{-x\MMa}] \; =\; \frac{1}{\Ga(\a)}\;  \EE[e^{-x\cMMa}]$$
where $\cMMa \elaw \cZa^{-\a}\!$ and $\cZa$ is the so-called size bias of order $-\a$ of $\Za$ that is the random variable with density function $\a\Ga(\a)\, x^{-\a}f^{}_{\Za}(x),$ which is characterized by
$$\EE[f(\cZa)] \; =\; \frac{\EE[\Za^{-\a} f(\Za)]}{\EE[\Za^{-\a}]}$$
for every $f$ bounded continuous. Finally, for every $\bt > \a,$ Lemma 2 in \cite{Mi} and a change of variable show that
$$\Eab(-x) \; =\;  \frac{1}{\Ga(\bt)}\;  \EE[e^{-x(\B_{\a,\bt-\a}^{\a}\times\,\cMMa)}].$$
In particular we observe the identity 
$$\MMa\, \elaw\,\B_{\a,1-\a}^{\a}\,\times\,\cMMa$$ 
for every $\a\in (0,1),$ which is not obvious at first sight. Notice also that the above identifications can be performed directly from the proof of Part (c) of the Theorem in \cite{Sch}.

\subsubsection{An improvement of {\em (a)}} It is a well-known fact - see e.g. Theorem 2 in \cite{MS} - that if $f$ is CM, then the function $f(x^\ga)$ is also CM for every $\ga\in (0,1).$ In particular, we see from Part (b) that the function $\cDab(x^{1/\a})$ is CM for every $\a\in (1, 2]$ and $\bt \ge 1.$ In this paragraph, we improve Part (a) and show that the function $\Dab(x^{1/\a})$ is also CM for every $\a\in (0,1)$ and $\bt\ge \a\vee (1-\a).$ The argument relies on the following lemma, which is fairly obvious. If $f$ is a positive function and $\mu$ a positive measure on $(0,+\infty),$ we set
$$f\,\odot\, \mu (t)\; =\; \int_0^\infty f(\frac{t}{s})\, \frac{\mu(ds)}{s}\cdot$$

\begin{LEM}
\label{Zaa}
Let $\a\in (0,1)$ and 
$$F(x)\; =\; \int_0^\infty e^{-xt} f(t)\, dt$$
be a {\em CM} function. Then $F(x^{1/\a})$ is also {\em CM} if and only if the Bernstein density factorizes into $f = f^{}_{\Za}\odot \mu$ for some positive $\sigma-$finite measure $\mu.$ 
\end{LEM}

\proof Suppose that the factorization holds. Two changes of variable and Fubini's theorem show that
\begin{eqnarray*}
F(x)\; =\;\int_0^\infty e^{-xt} \,f^{}_{\Za}\!\odot \mu (t)\, dt & = & \int_0^\infty\lpa\int_0^\infty  e^{-xts} f^{}_{\Za}(t) \,dt\rpa \mu(ds)\\
& = & \int_0^\infty e^{-x^\a s^\a} \mu(ds) \; =\; \int_0^\infty e^{-x^\a t} \,\mu^{[\a]}(dt),
\end{eqnarray*}
where $\mu^{[\a]}=x^\a[\mu]$ is a push-forward of $\mu$ and also a positive $\sigma-$finite measure. This shows that $F(x^{1/\a})$ is CM with Bernstein measure $\mu^{[\a]}.$ 

Suppose now that $F(x^{1/\a})$ is CM and let $\mu$ be its Bernstein measure. Setting $\mu^{1/\a}=x^{1/\a}[\mu],$ a comparison and two changes of variable show that 
\begin{eqnarray*}
\int_0^\infty e^{-x t} f(t)\, dt\; =\; \int_0^\infty e^{-x^\a t} \mu(dt)& = & \int_0^\infty\lpa\int_0^\infty  e^{-xt^{1/\a}s} f^{}_{\Za}(s) \,ds\rpa \mu(dt)\\
& = & \int_0^\infty\lpa\int_0^\infty  e^{-x t s} f^{}_{\Za}(s) \,ds\rpa \mu^{[1/\a]}(dt)\\
& = & \int_0^\infty e^{-x s} \,f^{}_{\Za}\odot\mu^{[1/\a]}(s)\, ds,
\end{eqnarray*}
with the same notation as above for the push-forward. We hence obtain the required factorization $f = f^{}_{\Za}\odot\mu^{[1/\a]},$ by uniqueness of the Laplace transform.

\endproof

We can now finish the proof. Suppose first that $\a > 1/2.$ From (\ref{12+}) we see that it is enough to show that $D_{\a,\a}(x^{1/\a})$ is CM.  From (\ref{daa}), we know that $f_{\Za}^{}$ factorizes the Bernstein density $d_{\a,\a}$ of $D_{\a,\a},$ and we can conclude by Lemma \ref{Zaa}. Assume next that $\a \le 1/2.$ The decomposition
$$\Dab(x^{1/\a}) \; = \; \frac{1}{\Ga(\a+\bt -1)}\int_0^1 (1-t)^{\a+\bt -2} t^{-\a} {\widetilde D}_{\a, 1-\a} (x^{1/\a} t) \, dt\; + \;\sum^{\na-1}_{k=1} \frac{x^{k-1/\a}}{\Ga(\a k+\bt -1)},$$
and the preceding argument entail that it is enough to show that $f_{\Za}^{}$ factorizes the Bernstein density of ${\widetilde D_{\a,1-\a}}.$ By (\ref{da1a}), the latter  is a constant multiple of 
$$\tha\; =\; f^{}_{\Za}\;\odot\;\tga\;\odot\; f^{}_{\B_{\a, 1-2\a}}\;\odot\;f^{}_{\Xa}\;\odot\;f^{}_{\G^{1/\a}_{1/\a}}\;\odot\;\lpa \bigodot_{k=1}^{\na -1} f^{}_{\U^{1/(1-\a k)}}\rpa,$$
and hence satisfies the required property.

\qed

\section{Proof of Theorem B}

Observe first that $E_{1, 1} (x) =  e^x = E_{1,1}' (x),$ so that the assertions (a) and (c) are obvious for $\a = 1.$ To handle the case $\a = 1, \bt > 1$ we appeal to the formul\ae
\begin{equation}
\label{Fab}
\Fab(x)\; =\; \frac{1}{\Ga(\bt -1)}\int_0^1 (1-t)^{\bt -2} F_{\a, 1} (xt) \, dt
\end{equation}
and
\begin{equation}
\label{FabL}
\Lb\Fab(x)\; =\; \frac{1}{\Ga(\bt -1)}\int_0^1 (1-t)^{\bt -2} F_{\a, 1}' (xt) \, dt
\end{equation}
which can be obtained as for (\ref{Frac1}). Setting $\a =1$ entails
\begin{equation}
\label{gainc}
E_{1,\bt} (x)\; =\;\Lb E_{1,\bt} (x)\; =\;\frac{1}{\Ga(\bt -1)}\int_0^1 (1-t)^{\bt -2} e^{xt} \, dt  \; =\; \frac{x^{1-\bt} e^x\ga(\bt-1,x)}{\Ga(\bt -1)},
\end{equation}
where the last equality follows from a straightforward change of variable. This shows that all functions in (a)-(d) are zero if $\a = 1.$ We now focus on the case $\a \neq 1.$ 

\subsection{Proofs of (a) and (b)} We first notice that (b) is a simple consequence of (a). Indeed, (\ref{Fab}), (\ref{FabL}) and the last equality in (\ref{gainc}) show that
$$\frac{x^{1-\bt} e^x\ga(\bt-1,x)}{\a\Ga(\bt -1)}\, -\,\Fab(x)\; =\; \frac{1}{\Ga(\bt -1)}\int_0^1 (1-t)^{\bt -2}\lpa \frac{e^{xt}}{\a} - F_{\a, 1} (xt)\rpa  dt$$
and
$$\Lb\Fab(x)\, -\,\frac{x^{1-\bt} e^x\ga(\bt-1,x)}{\a\Ga(\bt -1)}\; =\; \frac{1}{\Ga(\bt -1)}\int_0^1 (1-t)^{\bt -2}\lpa  F_{\a, 1}' (xt) - \frac{e^{xt}}{\a}\rpa  dt,$$
and it is clear that (a) entails that both functions between brackets are CM in $x$ for all $t \in [0,1].$ We hence focus on the case $\bt = 1.$ We will proceed again via Laplace inversion. Computing the Laplace transforms yields the following identity, which makes sense for every $s >0:$
$$\int_0^\infty e^{-sx} \lpa  \frac{e^{x}}{\a} - \Ea(x^\a)\rpa  dx\; =\; \frac{1}{\a (s -1)}\, -\, \frac{s^{\a-1}}{s^\a -1}\cdot$$
However, to show as in the preceding section that the function $\Fa$ on the right-hand side satisfies $\Re (e^{\i\theta/2} \Fa(re^{\i\theta})) \ge 0$ for every $r > 0$ and $\theta \in (-\pi,\pi),$ is very tedious. We hence follow a direct approach and compute the Stieltjes transform
$$\int_0^\infty \frac{dt}{s+t} \lpa \frac{\sin(\pi\a) t^{\a-1}}{\pi(t^{2\a} - 2\cos(\pi\a)t^\a +1)}\rpa  \;= \; \frac{\sin(\pi\a)}{\pi\a}\int_0^\infty \frac{dt}{(s+t^{1/\a})(t - e^{\i\pi\a})(t-e^{-\i\pi\a})},$$
which rewrites
$$\frac{1}{2\pi\i\a}\lpa\int_0^\infty \frac{dt}{(s+t^{1/\a})(t-e^{\i\pi\a})}\, -\,\int_0^\infty \frac{dt}{(s+t^{1/\a})(t-e^{-\i\pi\a})}\rpa.$$
Consider the contour $\Ga_R$ made out of the segment $[Re^{\i\pi(1-\a)}, Re^{-\i\pi\a}]$ oriented downwards, and of the half-circle $\cC_R$ leading anticlockwise from $Re^{-\i\pi\a}$ to $Re^{\i\pi(1-\a)}.$ Taking $s\neq 1, R > s,$ and applying the residue theorem shows on the one hand that
$$\int_{\Ga_R} \frac{dz}{(s-z^{1/\a})(z-1)}\; =\; 2\pi\i \lpa\frac{1}{s -1}\, -\, \frac{\a s^{\a-1}}{s^\a -1}\rpa.$$
On the other hand, the integral on the left-hand side is evaluated as
$$\int_0^R \frac{dt}{(s+t^{1/\a})(t-e^{\i\pi\a})}\; -\;\int_0^R \frac{dt}{(s+t^{1/\a})(t-e^{-\i\pi\a})}\; +\; \int_{\cC_R} \frac{dz}{(s-z^{1/\a})(z-1)}\cdot$$
Letting $R\to +\infty$ and putting everything together entails
$$\int_0^\infty \frac{dt}{s+t} \lpa \frac{\sin(\pi\a) t^{\a-1}}{\pi(t^{2\a} - 2\cos(\pi\a)t^\a +1)}\rpa  \;= \;\frac{1}{s -1}\, -\, \frac{\a s^{\a-1}}{s^\a -1}\cdot$$
Therefore, by Laplace inversion, 
\begin{equation}
\label{Exa}
\frac{e^x}{\a} - \Ea (x^\a)\; =\; \int_0^\infty e^{-xt} \lpa \frac{\sin(\pi\a) t^{\a-1}}{\pi(t^{2\a} - 2\cos(\pi\a)t^\a +1)}\rpa dt
\end{equation}
is a CM function as required.

\qed

\begin{REMS} {\em (a) The above formula (\ref{Exa}) should be compared with the classical formula (see (7.7) in \cite{HMS}, bewaring the misprint therein)
$$\Ea(-x^\a)\; =\; \int_0^\infty e^{-xt} \lpa \frac{\sin(\pi\a) t^{\a-1}}{\pi(t^{2\a} +2\cos(\pi\a)t^\a +1)}\rpa dt.$$
In particular setting $\a = 1/2$ we retrieve the identity $E_{1/2}(-\sqrt{x}) + E_{1/2}(\sqrt{x}) = 2e^x,$ which follows directly from the very definition of $E_{1/2}.$ \\
 
(b) Subtracting (\ref{Exa}) on both sides of (\ref{b2}), we get the other CM function
$$\a x^{\a -1}\Ea'(x^\a)\; -\; \frac{e^x}{\a}  \; =\;\int_0^\infty e^{-xt} \lpa \frac{\sin(\pi\a) t^\a}{\pi(t^{2\a} - 2\cos(\pi\a)t^\a +1)}\rpa dt,$$
which can also be obtained from a mere differentiation of (\ref{Exa}).}
\end{REMS}

\subsection{Proofs of (c) and (d)} Again, a direct application of (\ref{Fab}), (\ref{FabL}) and the last equality in (\ref{gainc}) show that (d) is a consequence of (c). To show (c), it is possible to use a contour integral analogous to that of (a). However, we will provide yet another argument which is specific to the case $\a\in (1,2].$ This proof is slightly lengthier but it has an independent interest. Referring to \cite{TS1} for details and further references, let $\{X_t, \, t\ge 0\}$ be the  spectrally positive L\'evy $\a-$stable process, starting from zero and normalized such that
\begin{equation}
\label{normal}
\EE\lcr e^{-s X_t}\rcr \; =\; e^{t s^\a}, \quad s,t \ge 0.
\end{equation}
Taking the Laplace transform, for every $\lbd > s^\a$ we get
$$\frac{1}{\lbd - s^\a}\; =\; \int_0^\infty e^{-\lbd t} \EE\lcr e^{-s X_t}\rcr dt\; =\; \int_\rl e^{-sx}\lpa \int_0^\infty  e^{-\lbd t} f^{}_{X_1} (xt^{-1/\a}) t^{-1/\a} dt\rpa dx,$$
where the second equality follows from Fubini's theorem and the $(1/\a)-$self-similarity of $\{X_t, \, t\ge 0\}.$ Differentiating with respect to $s$ and integrating in $\lbd$ we get
$$\frac{\a s^{\a-1}}{s^\a-\lbd}\; =\; \int_\rl e^{-sx}\lpa \int_0^\infty  e^{-\lbd t} f^{}_{X_1} (xt^{-1/\a}) xt^{-(1+1/\a)} dt\rpa dx.$$
Specifying to $\lbd =1$ and adding the two equalities, we obtain
$$\frac{s^{\a-1} -1}{s^\a- 1}\; =\; \int_\rl e^{-sx}\lpa \int_0^\infty  e^{-t} f^{}_{X_1} (xt^{-1/\a}) t^{-1/\a} (1+ \frac{x}{\a t})dt\rpa dx,$$
which makes sense for every $s > 0$. It is well-known - see Theorem 2.10.2 in \cite{Z} - that $\vert X_1\vert$ conditioned on $\{X_1 < 0\}$ has the same law as $\Za^{-\a},$ and that $\pb [X_1 < 0] = 1/\a.$ Hence, for every $x < 0$ one has
\begin{eqnarray*}
\int_0^\infty  e^{-t} f^{}_{X_1} (xt^{-1/\a}) t^{-1/\a} \, dt & = & \int_0^\infty  e^{-t} f^{}_{\Za} (t \vert x\vert^{-\a}) \frac{t\, dt}{\vert x\vert^{\a +1}}\\
& = & \vert x\vert^{\a -1}\int_0^\infty  e^{-\vert x\vert^\a t} f^{}_{\Za} (t) t \, dt \; =\; \frac{e^{-\vert x\vert}}{\a}\cdot
\end{eqnarray*}
Similarly, for every $x < 0$ we have 
\begin{eqnarray*}
-\lpa\int_0^\infty  e^{-t} f^{}_{X_1} (xt^{-1/\a}) \frac{x\, dt}{\a t^{1+1/\a}} \rpa & = & \frac{1}{\a \vert x\vert^\a}\int_0^\infty  e^{-t} f^{}_{\Za} (t \vert x\vert^{-\a}) dt\\
& = & \frac{1}{\a}\int_0^\infty  e^{-\vert x\vert^\a t} f^{}_{\Za} (t) dt \; =\; \frac{e^{-\vert x\vert}}{\a}\cdot
\end{eqnarray*}
Therefore, for every $s\ge 0,$
$$\frac{s^{\a-1} -1}{s^\a- 1}\; =\; \int_0^\infty e^{-sx}\lpa \int_0^\infty  e^{-t} f^{}_{X_1} (xt^{-1/\a}) t^{-1/\a} (1+ \frac{x}{\a t})dt\rpa dx, \qquad s\ge 0.$$
Inverting the Laplace transform, we deduce
\begin{eqnarray*}
\overline{D_{\a, 1}}(x)\; =\;F_{\a, 1}(x) \, -\, F_{\a,1}'(x) & = &\int_0^\infty  e^{-t} f^{}_{X_1} (xt^{-1/\a}) t^{-1/\a} (1+ \frac{x}{\a t})dt\\
& = & \int_0^\infty  e^{-x^\a t} f^{}_{X_1} (t^{-1/\a}) (1+ \a x^{\a -1} t)\frac{dt}{\a t^{1+1/\a}}\\
&=& G_\a (x)\; -\; G_\a'(x)
\end{eqnarray*}
where we have set $X_1^+ = X_1\, \vert\, X_1 > 0$ and 
$$G_\a(x)\; =\; \int_0^\infty  e^{-x^\a t} f^{}_{X_1} (t^{-1/\a})\frac{dt}{\a t^{1+1/\a}}\; =\; \lpa 1- \frac{1}{\a}\rpa \EE[ e^{- x^\a (X_1^+)^{-\a}}]$$
for every $x \ge 0.$ Solving the linear ODE with initial condition $F_{\a, 1} (0) = 1$ shows that
$$F_{\a, 1}(x)\; =\; \lpa 1- \frac{1}{\a}\rpa \EE[ e^{- x^\a (X_1^+)^{-\a}}]\; +\; \frac{e^x}{\a}\cdot$$
On the other hand, formula (3.3.16) in \cite{Z}  (beware the notation for $Z(\a,\rho)$ which is that of Chapter 3.1 therein) and Bochner's subordination for stable subordinators entail
$$(X_1^+)^{-\a}\; \elaw\; \lpa \frac{\Z_{1-1/\a}}{\Z_{\a -1}}\rpa^{\a -1}\; \elaw\; \lpa \frac{\Z_{\a -1}}{\Z_{\a -1}}\rpa^{\a -1}\!\times\; \Z_{1/\a}.$$
By Lemma \ref{Zaa}, this entails
$$\Ea(x^\a)\; -\; \frac{e^x}{\a}\; =\; \lpa 1- \frac{1}{\a}\rpa \EE[ e^{- x\U_{\a -1}^{1/\a}}]$$
with the notation
$$\U_{\a-1}\; =\; \lpa \frac{\Z_{\a -1}}{\Z_{\a -1}}\rpa^{\a-1}.$$
It is well-known - see e.g. (3.3.16) in \cite{Z} and compare with the positive part of a drifted Cauchy random variable - that $\U_{\a -1}^{1/\a}$ has the explicit density
$$\frac{\a \sin (\pi(\a -1)) t^{\a -1}}{\pi (\a-1) (t^{2\a} + 2 \cos (\pi(\a-1)) t^\a + 1)}\; =\; \frac{-\a \sin (\pi\a) t^{\a -1}}{\pi (\a-1) (t^{2\a} - 2 \cos (\pi\a) t^\a + 1)}\cdot$$
We finally obtain
\begin{equation}
\label{Eaxx}
\Ea(x^\a)\; -\; \frac{e^x}{\a}\; =\;\int_0^\infty e^{-xt} \lpa \frac{-\sin(\pi\a) t^{\a -1}}{\pi(t^{2\a} - 2\cos(\pi\a)t^\a +1)}\rpa dt,
\end{equation}
which concludes the proof. 

\qed

\subsection{List of the Bernstein measures and an improvement of (a) and (b)} We detail here the Bernstein densities associated with the CM functions of Theorem B, discarding the trivial case $\a = 1$ where all functions are zero. We will use all our previous notations without further repetition. We also improve parts (a) and (b) in the same way as we did in Paragraph 2.3.4. We next rephrase our results in the realm of spectrally positive stable L\'evy processes. Finally we compare some of them with the existing works on the so-called Mittag-Leffler distributions.

\subsubsection{The case $\a\in (0,1)$} From (\ref{Exa}) and (\ref{Fab}) we have
\begin{equation}
\label{Exx}
\frac{e^x}{\a}\; -\; \Ea(x^\a)\; =\; \lpa\frac{1}{\a} - 1\rpa \EE[e^{-x \U_{1-\a}^{1/\a}}]
\end{equation}
and, for every $\bt > 1,$
$$\frac{x^{1-\bt} e^x\ga(\bt-1,x)}{\a\Ga(\bt -1)}\, -\,\Fab(x)\; =\; \frac{1}{\Ga(\bt)}\lpa\frac{1}{\a} - 1\rpa \EE[e^{-x (\U_{1-\a}^{1/\a}\times\, \B_{1,\bt-1})}].$$
From (\ref{dacis}), we finally see that the infinite Bernstein measure associated with the CM function
$$\Lb\Fab(x)\, -\,\frac{x^{1-\bt} e^x\ga(\bt-1,x)}{\a\Ga(\bt -1)}$$
has density 
$$\frac{1}{\Ga(\bt)}\, f_ {\B_{1,\bt-1}}^{}\,\odot\, {\hat \fa},$$
where 
$${\hat \fa} (t)\; =\; \frac{\sin(\pi\a) t^\a}{\pi(t^{2\a} - 2\cos(\pi\a)t^\a +1)}\cdot$$

\subsubsection{The case $\a\in (1,2]$} We have seen above that
$$\Ea(x^\a)\; -\; \frac{e^x}{\a} \; =\; \lpa 1- \frac{1}{\a}\rpa \EE[e^{-x \U_{\a-1}^{1/\a}}],$$
and again (\ref{Fab}) entails
$$\Fab (x)\, -\, \frac{x^{1-\bt} e^x\ga(\bt-1,x)}{\a\Ga(\bt -1)}\; =\; \frac{1}{\Ga(\bt)}\lpa 1-\frac{1}{\a}\rpa \EE[e^{-x (\U_{\a-1}^{1/\a}\times\, \B_{1,\bt-1})}]$$
for every $\bt > 1.$ Introduce now with the size bias of order 1 of the random variable 
$\U_{\a-1}^{1/\a}$ and denote it by $(\U_{\a-1}^{1/\a})^{(1)}.$ This is a proper random variable with density
$$\frac{-\a\sin(\pi\a) t^\a}{\pi(t^{2\a} - 2\cos(\pi\a)t^\a +1)}\cdot$$
The above entail the further formul\ae
$$\frac{e^x}{\a} \; -\; \a x^{\a -1}\Ea'(x^\a)\; =\; \frac{1}{\a}\, \EE[e^{-x (\U_{\a-1}^{1/\a})^{(1)}}]$$
and, for every $\bt > 1,$
$$\frac{x^{1-\bt} e^x\ga(\bt-1,x)}{\a\Ga(\bt -1)}\, -\,\Lb\Fab (x)\; =\; \frac{1}{\a\Ga(\bt)}\,\EE[e^{-x ((\U_{\a-1}^{1/\a})^{(1)}\times\, \B_{1,\bt-1})}].$$

\subsubsection{An improvement of {\em (a)} and {\em (b)}} In this paragraph we show the following proposition.

\begin{PROP} 
\label{prp}
Let $\a\in (0,1).$ The function
$$x\;\mapsto\;\frac{e^{x^{1/\a}}}{\a}\, -\,\Ea(x)$$ 
is {\em CM} if and only if $\a \ge 1/2.$ 
\end{PROP}

\proof We first compute the fractional moments of the random variable $\U_{1-\a}^{1/\a}$ appearing in (\ref{Exx}). For every $s\in (-\a, \a),$ we find
\begin{eqnarray*}
\EE[ \U_{1-\a}^{s/\a}] & = & \frac{\Ga (1-\frac{s}{\a})\Ga(1+\frac{s}{\a})}{\Ga (1 - (\frac{1-\a}{\a})s)\Ga(1+(\frac{1-\a}{\a})s)}\\
& = & \EE[\Za^s]\;\times\;\lpa \frac{\Ga (1-s)}{\Ga (1 - (\frac{1-\a}{\a})s)}\;\times\;\frac{\Ga(1+\frac{s}{\a})}{\Ga(1+(\frac{1-\a}{\a})s)}\rpa.
\end{eqnarray*}
If $\a \ge 1/2,$ this can be read off in the following way
$$\EE[ \U_{1-\a}^{s/\a}] \; =\;  \EE[\Za^s]\;\times\; \EE[\Wa^{\frac{s}{\a}}],$$
with the notation
$$\Wa \; =\; \lpa \frac{\Z_{\frac{1-\a}{\a}}}{\Z_{1-\a}}\rpa^{1-\a}.$$
From (\ref{Exx}) and Lemma \ref{Zaa}, this entails that
$$ \frac{e^{x^{1/\a}}}{\a}\, -\,\Ea(x)\; =\; \lpa \frac{1}{\a} - 1\rpa \EE [e^{-x {\bf W}_\a}]$$ 
is a CM function. On the other hand, if $\a < 1/2$ then the quotient under brackets vanishes at $s = \a/(1-\a) < 1.$ If this quotient were the Mellin transform of a positive measure, then it would vanish inside its definition strip which is $-\a < \Re(s) < 1,$ and this is a contradiction. This shows that $f_{\Za}$ does not factorize the Bernstein density $f_{\U_{1-\a}^{1/\a}}$ and, by Lemma \ref{Zaa}, that 
$$x\; \mapsto\;\frac{e^{x^{1/\a}}}{\a}\, -\,\Ea(x)$$ 
is not CM.

\endproof

From (\ref{Fab}) and the above, we see that for every $\a \in [1/2,1)$ the function
$$\frac{x^{\frac{1-\bt}{\a}} e^{x^{\frac{1}{\a}}}\ga(\bt-1,x^{\frac{1}{\a}})}{\a\Ga(\bt -1)}\; -\;\Eab (x) \; =\; \frac{1}{\Ga(\bt)}\lpa 1-\frac{1}{\a}\rpa \EE[e^{-x (\Wa\,\times\, \B_{1,\bt-1})}]$$
is CM, which also improves (b). Observe also that this function taken at $x^\ga$ is not CM for any $\ga > 1,$ since the derivative at zero is then zero. It is however unclear whether this function is CM or not for $\a < 1/2.$ To characterize the CM property of the functions
$$\a x^{\a -1}\Ea'(x)\; -\; \frac{e^{x^{1/\a}}}{\a}\qquad\mbox{and}\qquad
\Lb\Fab (x^{\frac{1}{\a}})\; -\; \frac{x^{\frac{1-\bt}{\a}} e^{x^{\frac{1}{\a}}}\ga(\bt-1,x^{\frac{1}{\a}})}{\a\Ga(\bt -1)}$$
is also an open question.

\subsubsection{Connections with spectrally positive stable L\'evy processes} In this paragraph we fix $\bt = 1$ and $\a\in (1,2].$ As in Section 3.2, let $X = \{X_t, \, t\ge 0\}$ be the  spectrally positive L\'evy $\a-$stable process, starting from zero and normalized by (\ref{normal}). Let
$$T_1\; =\; \inf\{t > 0,\; X_t \, > \, 1\}\qquad\mbox{and}\qquad \tau_1\; =\; \inf\{t > 0,\; X_t \, = \, 1\}$$
be respectively the first passage time above 1 and the first hitting time of 1 for $X.$ It is well-known from the Wiener-Hopf factorization - see (7) in \cite{TS0} and the references therein for further details - that 
$$\Ea(x) \; -\; \a x^{1-1/\a} \Ea'(x)\; =\;\EE[e^{-x T_1}].$$
On the other hand, a consequence of Fristedt's formula  - see (1.4) in \cite{TS1} and the computations therebefore for an explanation - is that 
$$\frac{e^{x^{1/\a}}}{\a} \; -\; \a x^{1-1/\a} \Ea'(x)\; =\;\frac{1}{\a}\;\EE[e^{-x \tau_1}].$$
Besides, we saw above in Paragraph 3.3.1 that
\begin{eqnarray*}
\Ea(x)\; -\; \frac{e^{x^{1/\a}}}{\a} & = & \lpa 1- \frac{1}{\a}\rpa \EE[e^{-(x \U_{\a-1})^{1/\a}}]\\
& = & \lpa 1- \frac{1}{\a}\rpa \EE[e^{-x (\U_{\a-1}\times\Z_{1/\a})}]\; =\; \lpa 1- \frac{1}{\a}\rpa \EE[e^{-x (X^+_1)^{-\a}}],
\end{eqnarray*}
where the third equality follows from formula (3.3.16) in \cite{Z} and Bochner's subordination. Hence, we see that the decomposition in two CM functions 
$$\Ea(x) \; -\; \a x^{1-1/\a} \Ea'(x)\; =\; \lpa \Ea(x)\; -\; \frac{e^{x^{1/\a}}}{\a}\rpa\; +\; \lpa\frac{e^{x^{1/\a}}}{\a} \; -\; \a x^{1-1/\a} \Ea'(x)\rpa$$
can be interpreted as an explicit relationship between three Laplace transforms connected to the L\'evy process $X$:
$$\EE[e^{-x T_1}]\; =\; \lpa 1- \frac{1}{\a}\rpa \EE[e^{-x (X^+_1)^{-\a}}]\; +\; \frac{1}{\a}\;\EE[e^{-x \tau_1}].$$
In particular, inverting these Laplace transforms shows that for every $x\ge 0,$
\begin{equation}
\label{Int}
f^{}_{T_1} (x)\; =\; \lpa 1- \frac{1}{\a}\rpa f^{}_{(X^+_1)^{-\a}} (x)\; +\; \frac{1}{\a}\;f^{}_{\tau_1} (x),
\end{equation}
which clarifies formula (14) in \cite{TS0}. Notice that the decomposition (\ref{Int}) can actually also be derived in comparing Propositions 2 and 3 in \cite{TS1}, Corollary 6 in \cite{TS1} and the classical series representation (2.4.6) in \cite{Z}.

\begin{REMS} {\em Set now $\bt = 1$ and $\a \in [1/2,1).$ Let $Y = \{Y_t, \, t\ge 0\}$ be the  spectrally positive L\'evy $(1/\a)-$stable process, starting from zero and normalized as in (\ref{normal}). Recall that this process has positivity parameter $\pb[Y_1 > 0] = 1-\a.$ It follows from the proof of Proposition \ref{prp}, and the same considerations as above around formula (3.3.16) in \cite{Z}, that
\begin{equation}
\label{Last}
\frac{e^{x^{1/\a}}}{\a}\, -\,\Ea(x)\; =\; \lpa \frac{1}{\a} - 1\rpa \EE [e^{-x {\bf W}_\a}]\; =\; \lpa \frac{1}{\a} - 1\rpa \EE [e^{-x Y_1^+}].
\end{equation}
This hence establishes a link between $\Ea$ and yet another spectrally positive stable L\'evy process. The latter is actually a consequence of the usual connection between Mittag-Leffler functions and spectrally negative stable L\'evy processes. Consider indeed $I_1 = \inf\{ Y_t, \, t\in [0,1]\}.$ It is well-known - see e.g. Proposition 1 (iii) in \cite{Bi} - and has been used during the proof of Theorem B (c) that 
$$\Ea(x)\; =\; \EE[e^{x I_1}]\; =\; \EE[e^{x Y_1^-}],\qquad x\in \rl,$$
with $Y_1^- =\vert Y_1\vert$ conditioned on $\{Y_1<0\}.$ The latter identity entails (\ref{Last}) with the help of (\ref{normal}) and of the basic identity
$$\EE[e^{-x Y_1}]\; =\; \a \EE[e^{x Y_1^-}] \; +\; (1-\a)\EE [e^{-x Y_1^+}], \qquad x\ge 0.$$ }
\end{REMS} 

\subsubsection{Comparison with the other Mittag-Leffler distribution} Let $\L$ be the unit exponential random variable and $\a\in (0,1)$. The random variable 
\begin{equation}
\label{Mla}
\Mla\; =\; \L^{1/\a}\;\times\; \Za
\end{equation}
is known as the "other" Mittag-Leffler random variable, the classical one being the above $I_1.$ It is a particular instance of the Linnik (or geometric stable) random variables. It has an explicit and completely monotone density which is 
$$x^{\a -1}E_{\a,\a}(-x^\a),$$ 
and an explicit Laplace transform 
$$\EE [e^{-\lbd \Mla}]\; =\; \frac{1}{1+\lbd^\a}\cdot$$
It is infinitely divisible viz. the function
$$\lpa \frac{1}{1+\lbd^\a}\rpa^t$$
is CM for all $t\ge 0,$ and the associated semi-group, which is sometimes called the Mittag-Leffler semi-group in the literature, has also a semi-explicit transition density which is given by (20.1.2) in \cite{HMS}. The factorization (\ref{Mla}) means that this semi-group is subordinated to the positive $\a-$stable semi-group. We refer to Section 19 in \cite{HMS} for details and references on the above properties, and also for further features of the other Mittag-Leffler distribution. \\

Let now $\a\in (1,2)$. In the present paper we exhibited a positive random variable $\CMla,$ with an explicit and completely monotone density which is 
$$\Ea(x^\a)\; -\; \a x^{\a-1}\Ea'(x^\a),$$ 
and an explicit Laplace transform 
$$\EE [e^{-\lbd \CMla}]\; =\; \frac{\lbd^{\a -1} -1}{\lbd^\a -1}\cdot$$
This random variable is infinitely divisible and the associated semi-group is  subordinated to the positive $(\a/2)-$stable semi-group, viz. one has the factorization
$$\CMla\; =\; \X^{2/\a}\;\times\; \Z_{\a/2},$$
where the random variable $\X$ is infinitely divisible. One can also show that $\CMla$ has ${\bf ML}_{\a/2}$ as an additive factor. In a forthcoming work, we will show these properties and also present further features of this Mittag-Leffler random variable of the second kind.


\begin{thebibliography}{10}


\bibitem{Bi}
N.~H.~Bingham. Maxima of sums of random variables and suprema of stable
processes. {\em Z. Wahrsch. verw. Gebiete} {\bf 26}, 273-296, 1973. 

\bibitem{D}
M.~M.~Djrbashian. {\em Integral transforms and representations of functions in the complex plane.} Nauka, Moskva, 1966.

\bibitem{E3} 
A.~Erdelyi. {\em Higher transcendental functions. Vol III.} McGraw-Hill, New-York, 1953.

\bibitem{HMS} 
H.~J. Haubold, A.~M.~Mathai and R.~K.~Saxena. Mittag-Leffler functions and their applications. {\em J. Appl. Math.}, 51 pages, 2011.

\bibitem{HW}
I.~I.~Hirschman and D.~V.~Widder. {\em The convolution transform.} Princeton University Press, Princeton, 1955.

\bibitem{Mai}
F.~Mainardi. {\em Fractional calculus and waves in linear viscoelasticity.} Imperial College Press, London, 2010.

\bibitem{Mi}
K.~S.~Miller and S.~G.~Samko. A note on the complete monotonicity of the generalized Mittag-Leffler function. {\em Real Anal. Exch.} {\bf 23} (2), 753-756, 1999.
 
\bibitem{MS}
K.~S.~Miller and S.~G.~Samko. Completely monotonic functions. {\em Integr. Transf. Spec. Funct.} {\bf 12} (4), 389-402, 2001.

\bibitem{P}
H.~Pollard. The completely monotonic character of the Mittag-Leffler function $\Ea(-x)$. {\em Bull. Am. Math. Soc.} {\bf 52}, 908-910, 1948.

\bibitem{Sch}
W.~R.~Schneider. Completely monotone generalized Mittag-Leffler functions. {\em Expo. Math.} {\bf 14}, 3-16, 1996.

\bibitem{TS0}
T.~Simon. Fonctions de Mittag-Leffler et processus de L\'evy stables sans sauts n\'egatifs. {\em Expo. Math.} {\bf 28}, 290-298, 2010.

\bibitem{TS1}
T.~Simon. Hitting densities for spectrally positive stable processes. {\em Stochastics} {\bf 83} (2), 203-214, 2011.

\bibitem{W}
D.~V.~Widder. {\em The Laplace transform.} Princeton University Press, Princeton, 1946. 

\bibitem{WZ}
R.~Wong and Y.~Q.~Zhao. Exponential asymptotics of the Mittag-Leffler function. {\em Constr. Approx.} {\bf 18}, 355-385, 2002.

\bibitem{Z}
V.~M.~Zolotarev. {\em One-dimensional stable distributions.} Nauka, Moskva, 1983. 

\end{thebibliography}
\end{document}